\documentclass{amsart}
\usepackage{amsthm, amsmath}

    \newtheorem{theorem}{Theorem}[section]
    \newtheorem{corollary}[theorem]{Corollary}
    \newtheorem{proposition}[theorem]{Proposition}
    \newtheorem{remark}[theorem]{Remark}
    \newtheorem{lemma}[theorem]{Lemma}

\begin{document}
\markboth{D.S. Ramana \& O. Ramar\'e}
{A Variant of the Truncated Perron's Formula and Primitive Roots}

\title{A Variant of the Truncated Perron's Formula and Primitive Roots}

\author{D.S. Ramana}
\address{HBNI/Harish-Chandra Research Institute, Jhunsi,\\
Allahabad -211 019, India.
\\    suri@hri.res.in}

\author{O. Ramar\'e}
\address{CNRS/Institut de Math{\'e}matiques de Marseille,
        \\Aix Marseille Universit{\'e}, Centrale Marseille, 
        \\Site Sud, Campus de Luminy, Case 907\\
13288 Marseille Cedex 9, France.\\
        ramare@math.univ-lille1.fr}

\maketitle
\begin{abstract}
  We show under the Generalised Riemann Hypothesis that for every
  $\delta>0$, almost every prime $q$ in $[Q,2Q]$ has the expected of
  prime primitive roots in the interval $[x,x+x^{\frac{1}2+\delta}]$
  provided $Q$ is not more than $x^{\frac{2}{3}-\epsilon}$. We obtain
  this via a variant of the classical truncated Perron's formula for
  the partial sums of the coefficients of a Dirichlet series.
 \end{abstract}


\section{Introduction}
\label{intro}

The classical truncated Perron's formula relates, for any $x \geq 1$,
the partial sum $\sum_{1 \leq n \leq x} a_n$ of the coefficients of a
Dirichlet series $F(s) = \sum_{n \geq 1} \frac{a_n}{n^s}$ with a
finite abscissa of convergence $\sigma_{c}$ to the integral on the
line segment $[\kappa - iT, \kappa + iT]$ of
$\frac{F(s)x^s}{2\pi i s}$, for any $T >0$ and
$\kappa >\max(0, \sigma_{a})$, where $\sigma_{a}$ is the abscissa of
absolute convergence of $F(s)$. The difference between these two
quantities is estimated by an error term that depends on a sum of the
absolute values $|a_n|$ of the $a_n$. We present here a variant that
has sums of the $a_n$ rather than $|a_n|$ and is valid for
$\kappa > \max(0,\sigma_{c})$. The basic version of this variant is
stated in Theorem~\ref{pp1}. This proposition results from a
simple rewriting of the Fourier adjunction formula
\begin{equation}
\label{ff}
\int_{{\bf R}} f(u) \hat{\phi}(u)  du =
\int_{{\bf R}} \hat{f}(u)\phi(u) du,
\end{equation} 
valid for any $f, \phi$ in $\text{L}^{1}({\bf R})$, applied with
$f(u) = e^{-\kappa u} \sum_{1\leq n \leq xe^{u}} a_n $ and suitable~$\phi$. Corollaries~\ref{pp2} puts Theorem~\ref{pp1} in applicable form. These are stated with the aid of
notation introduced at the head of Section \ref{var}. At the end of
this section we include a brief comparative description with other
variants of the Perron formula in the literature such as those in
G. Coppola \& S. Salerno \cite{cs1}, \cite{cs2}, J. Kaczorowski \&
A. Perelli\cite{kp}, J. Liu \& Y. Ye \cite{ly}, and Wolke\cite{wol}.
As an illustration of our version of the truncated Perron's formula,
we shall obtain the following result in Section 3.

\begin{theorem}
\label{theo1}
For any integer real numbers $Q\geq 2$, $\delta\le 1/3$ and $\theta\le 1/4$  and assuming the Generalised Riemann Hypothesis we have that for all but
  $O(Q^{\frac{11}{12}})$ 
  primes $q$ in $[Q,2Q]$, the number of prime primitive roots modulo $q$ in
  $[x, x+x^{\frac12+\delta}]$ is asymptotic to
  \begin{equation*}
    \frac{\varphi(q-1)}{q-1}\frac{x^{\frac12+\delta}}{\log x}
  \end{equation*}
  provided that $x\ge Q^{\frac{3}2+\theta}$.
  Furthermore, for almost all
  prime moduli $q$ in $[Q,2Q]$, the sum $\sum_{p}\mu(p)$ where $p$ ranges the
  primes from $[x, x+x^{\frac12+\delta}]$ that are primitive roots modulo~$q$
  is $o(x^{\frac12+\delta}/\log x)$ provided again that $x\ge Q^{\frac{3}2+\theta}$. 
\end{theorem}

Note that the modulus $q$ may be larger than the size of the interval
$x^{\frac12+\delta}$. The restriction to prime $q$ is only for
simplicity. When we are only interested in existence rather than an
asymptotic, sieve techniques may be employed to obtain much better
results as, for instance, in G. Martin \cite{mar}, where a bound for
the least prime primitive root is given under the GRH.

In the final section of this note, Section \ref{mli}, we consider the
effect of ``moving the line of integration $\sigma = \kappa$" in the
integrals on this line on the right hand sides of the formulae
supplied by Corollary~\ref{pp2}. In the
classical case the kernel $\phi$ is identically equal to 1 on
$[\kappa - iT, \kappa + iT]$. Our choices for $\phi$ are, however,
sufficiently smooth, compactly supported, piecewise polynomial
functions on $[\kappa - iT, \kappa + iT]$. These functions extend
holomorphically in horizontal strips and, in general, these extensions
are incompatible on adjacent strips. Nevertheless, Proposition
\ref{pp3} tells us that the smoothness of $\phi$ is enough to
guarantee that the error due to this incompatibility on moving the
line of integration is $O({x^{\kappa}}/{T^2})$ under resonable
assumptions on $F$.

\vspace{2mm}
\noindent
Throughout this article we use $e(z)$ to denote $e^{2\pi i z}$, for
any complex number $z$. Further, all constants implied by the symbols
$\ll$ and $\gg$ are absolute except when dependencies are indicated,
either in words or by subscripts to these symbols. We will use the
terms majorised and minorised to mean $\gg$ and $\ll$
respectively. The Fourier transform $\widehat{f}$ of an integrable
function $f$ on ${\bf R}$ is defined by
$\widehat{f}(u) = \int_{\bf R} f(t)e(-ut) dt$.

\section{The Variant}
\label{var}

Throughout this section, we let
$F(s) = \sum_{n \geq 1} {a_n}/{n^{s}}$ be a Dirichlet series with
a finite abscissa of convergence ${\sigma_{c}}$ and an abscissa of
absolute convergence $\sigma_{a}$. Also, let
$\sigma_{0} = \max{(0, \sigma_{c})}$ and for any real
$\sigma > \sigma_0$, let
$B(\sigma) = \sup_{N \geq 1}|\sum_{1 \leq n \leq N}
\frac{a_n}{n^{\sigma}}|$.
Then on writing $a_n$ as $\frac{a_n}{n^{\sigma}}\cdot n^{\sigma}$ and
using the Abel summation formula we obtain the classical bound of
E. Cahen \cite{cah}:
\begin{equation}
\label{cah}
 \Bigl|\sum_{1 \leq n \leq x} a_n  \Bigr| \leq 2B(\sigma) x^{\sigma}  ,
\end{equation}
valid for all $x \geq 1$ and any $\sigma > \sigma_0$. The following
theorem uses a test function $\phi$ and its Fourier transform
$\hat{\phi}$ to express $\sum_{1 \leq n \leq x} a_n$ in terms of
$F(s)$.

\begin{theorem} \label{pp1} Let $\phi$ a function in
  $\text{L}^{1}({\bf R})$ with $\phi(0) = 1$ and such that
  $\hat{\phi}$ is also in $\text{L}^{1}({\bf R})$. Then for any
  $\kappa > \sigma_{0}$ and $x \geq 1$ we have
\begin{multline}
\label{prop1}
\sum_{1 \leq n \leq x} a_n \; =\; \frac{1}{2\pi i} \int_{\kappa-i\infty}^{\kappa + i\infty} F(s)\phi\left(\frac{s-\kappa}{2\pi i }\right) \frac{x^s}{s} \; ds 
\\
+\int_{{\bf R}} \left(\sum_{1 \leq n \leq x} a_n - e^{-\kappa u}\sum_{1 \leq n \leq xe^{u}} a_n \right) \hat{\phi}(u) \; du .
\end{multline}
\end{theorem}

\begin{proof}
For any $f$ in $\text{L}^{1}({\bf R})$ we have 
\begin{equation}
\label{eq1}
 f(0) =  \int_{\bf R} f(0) \hat{\phi}(u) du =  \int_{{\bf R}} \hat{f}(u) \phi(u) du + \int_{{\bf R}} (f(0)- f(u)) \hat{\phi}(u) du \; ,
\end{equation}
on taking account of $\int_{\bf R} \hat{\phi}(u) du = \phi(0) = 1$ and  \eqref{ff}, valid since $\phi$ and  $\hat{\phi}$ are also in $\text{L}^{1}({\bf R})$. The relation \eqref{prop1} results on using \eqref{eq1} with 
\begin{equation*}
f(u) = e^{-\kappa u} \sum_{1\leq n \leq xe^{u}} a_n \; ,
\end{equation*}
for the given $\kappa > \sigma_0$ and $x \geq 1$. Indeed, then $f(u) = 0$ for $u < -\log x$. Further, we have from \eqref{cah} with $\sigma$ in $(\sigma_{0}, \kappa)$ that $f$ is integrable in a neighbourhood of $+\infty$. Thus $f$ is in
$\text{L}^{1}({\bf R})$. A comparision of the right hand side of \eqref{prop1} with the last term of \eqref{eq1} now shows that it only remains to verify that 
\begin{equation}
\label{fou}
 \hat{f}(u) = \frac{x^{s}F(s)}{s} \;\;\;\;\text{where $s = \kappa + 2\pi i u$ ,} 
\end{equation}
for all $u \in {\bf R}$,
which is a well-known fact. For the sake of completeness, however, we provide a proof.
For any integer $m \geq 1$, let  
$f_{m}(u) = e^{-\kappa u} \sum_{1 \leq n \leq m} a_n \chi_{n}(u)$,  where $\chi_{n}(u)$ is 1 when $n \leq xe^{u}$ and is 0 otherwise. Then we 
certainly have 
$\lim_{m \rightarrow +\infty} f_{m}(u) = f(u)$ for all $u$ in ${\bf R}$. Also, (\ref{cah}) gives for any $\sigma$ in $(\sigma_{0}, \kappa)$ the bound 
\begin{equation*}
|f_{m}(u)| \leq 2B(\sigma)x^{\sigma} e^{(\sigma-\kappa)u}
\end{equation*}
for all $m \geq 1$ and all $u$ in ${\bf R}$. This allows us to apply the  dominated convergence theorem to justify the relation 
\begin{align}
\int_{{\bf R}} f(t) e^{-2\pi i u t} dt &= \notag
\lim_{m \rightarrow +\infty} \int_{{\bf R}} f_{m}(t) e^{-2\pi iu t} dt 
\\&=
\lim_{m \rightarrow +\infty} \sum_{1 \leq n \leq m} a_n \int_{{\bf R}} \chi_{n}(t) e^{-\kappa t} e^{-2\pi i ut} dt \; ,\label{eq3}
\end{align}
for all $u$ in ${\bf R}$. Since $\kappa > 0$ we have 
\begin{equation*}
\int_{{\bf R}} \chi_{n}(t) e^{-\kappa t} e^{-2\pi i ut} dt = \int_{ \log(\frac{n}{x})}  e^{-(\kappa + 2\pi i u) t} dt = \frac{x^{s}}{s\, n^{s}} \; , 
\end{equation*}
for all $n \geq 1$, where $s = \kappa + 2 \pi i u$. Also, since $\kappa > \sigma_{c}$ we have $\lim_{m \rightarrow +\infty} \sum_{1 \leq n \leq m} \frac{a_n}{n^s} = F(s)$. Consequently, \eqref{eq3} yields (\ref{fou}).
\end{proof}

The following corollary puts the second term on the right hand side
of \eqref{prop1} into a convenient form, with additional hypotheses on
the test function $\phi$. These hypotheses are satisfied when $\phi$
is a sufficiently smooth positive compactly supported function with
$\phi(0)=1$, as will be the case in our application.

\begin{corollary} \label{pp2} Let $\phi$ a function in $\text{L}^{1}({\bf R})$ with $\phi(0) = 1$ and such that 
  \begin{itemize}
  \item $\hat{\phi}$ is in $\text{L}^{1}({\bf R})$ and
 $\hat{\phi}(u) =  \hat{\phi}(-u)$ for all $u$ in ${\bf R}$.
  \item There is $m \geq 2$ such that  $C_{k}(\phi) = \displaystyle\sup_{u \in {\bf R}} |u^k \hat{\phi}(u)| < +\infty$ for  $0 \leq k \leq m +1$.
  \end{itemize}
  We set $C(\phi) = \max\limits_{0 \leq k \leq m+1} C_{k}(\phi)$. Then for
  any $\kappa > \sigma_{0}$, $x \geq 1$ and $T \ge1$ we have
 \begin{multline}
 \label{cor1}
 \sum_{1 \leq n \leq x} a_n \; = \frac{1}{2\pi i} \int_{\kappa-i\infty}^{\kappa + i\infty} F(s)\phi\left(\frac{s-\kappa}{2\pi iT }\right) \frac{x^s}{s} \; ds \, 
           \\+\int_{0}^{T} \left( \sum_{xe^{-\frac{u}{T}} < n \leq xe^{\frac{u}{T}}} a_n {\rm sgn}(x-n)\right) \hat{\phi}(u) \; du 
           \\ + {\mathcal O}^{*}\left( \frac{4C(\phi)(1+\kappa)^2e^{2\kappa} B(\kappa)x^{\kappa}}{T}\right).
 \end{multline}
\end{corollary}

\begin{proof}
  We first prove that
  \begin{multline}
 \label{prop2}
 \sum_{1 \leq n \leq x} a_n \; = \frac{1}{2\pi i} \int_{\kappa-i\infty}^{\kappa + i\infty} F(s)\phi\left(\frac{s-\kappa}{2\pi iT }\right) \frac{x^s}{s} \; ds \, 
            \\+\int_{0}^{T} \left(\sum_{xe^{-\frac{u}{T}} < n \leq xe^{\frac{u}{T}}} a_n {\rm sgn}(x-n)\right) \hat{\phi}(u) \; du 
            + \frac{\kappa}{T} \int_{0}^{T} \left(\sum_{xe^{-\frac{u}{T}} < n \leq xe^{\frac{u}{T}}} a_n\right) u\hat{\phi}(u) \; du 
            \\+ {\mathcal O}^{*}\left( \frac{4C(\phi)(1+\kappa^2)e^{2\kappa} B(\kappa)x^{\kappa}\log(eT)}{T^2}\right).
 \end{multline}
  To do so, we apply \eqref{prop1} to the
  function $u \mapsto \phi(\frac{u}{T})$, which we denote by
  $\psi$. Plainly, the first terms on the right hand sides of
  \eqref{prop2} and \eqref{prop1} applied to $\psi$ are the same. If
  for the given $\kappa > \sigma_{0}$ and $x \geq 1$ we set
  $A(u) = \sum_{1 \leq n \leq xe^u} a_n$ for all $u \in {\bf R}$ then,
  since $\hat{\psi}$ is also even by $(ii)$, the second term on the
  right hand side of \eqref{prop1} applied to $\psi$ can be written
 \begin{equation}
 \label{prop21}
 \int_{0}^{+\infty} \left(2A(0) - e^{-\kappa u}A(u) - e^{\kappa u}A(-u)\right)\hat{\psi}(u) du.  
 \end{equation}
First we estimate the contribution to the integral \eqref{prop21} from the interval $[1, +\infty)$. From \eqref{cah} with $\sigma = \kappa$ we see that $|A(0)|$ and $|e^{-\kappa u}A(u)|$ for $u \geq 1$ do not exceed $2B(\kappa)x^{\kappa}$. Since $A(-u) = 0$ when $u > \log x$, we similarly obtain $|e^{\kappa u}A(-u)| \leq 2B(\kappa)x^{\kappa}$ for $u \geq 1$. Consequently, we have 
 \begin{multline}
 \label{prop24}
 \biggl|\int_{1}^{+\infty} \left(2A(0) - e^{-\kappa u}A(u) - e^{\kappa
     u}A(-u)\right)\hat{\psi}(u) du\biggr|
 \\\leq 8B(\kappa)x^{\kappa}\int_{1}^{+\infty} |\hat{\psi}(u)| du \leq  4B(\kappa)x^{\kappa}C_3(\psi). 
 \end{multline}
Let us now define $h(z)$ for any complex number $z$ by $e^z = 1 +z +h(z)$. Then the contribution  to the integral \eqref{prop21} from the interval $(0,1)$ can be written as 
\begin{equation}
\label{prop22}
\begin{split}
 &\int_{0}^{1} \left(2A(0) - A(u) - A(-u)\right)\hat{\psi}(u) du 
 +\kappa \int_{0}^{1} \left(A(u)- A(-u)\right)u\hat{\psi}(u) du \\
 &-\int_{0}^{1} \left(h(-\kappa u)A(u) + h(\kappa u)A(-u)\right)\hat{\psi}(u) du . 
\end{split}
\end{equation}   
We estimate the third integral in \eqref{prop22} by means of the bounds $|h(z)| \leq \frac{|z|^2e^{|z|}}{2}$ for all $z \in {\bf C}$ and $|A(u)| \leq 2B(\kappa)x^{\kappa}e^{\kappa}$ for $|u| \leq 1$. The first of these bounds follows from the Taylor expansion of $e^z$ while the second follows from \eqref{cah} with $\sigma = \kappa$. We obtain 
\begin{multline}
\label{prop23}
\biggl|\int_{0}^{1} \left(h(-\kappa u)A(u) + h(\kappa
  u)A(-u)\right)\hat{\psi}(u) du\biggr|
\\ \leq
2\kappa^2 e^{2\kappa} B(\kappa)x^{\kappa}\int_{0}^{1}u^2|\hat{\psi}(u)| du.
\end{multline}
We have $\hat{\psi}(u) = T\hat{\phi}(uT)$ and therefore
$C_3(\psi) = \frac{C_3(\phi)}{T}$ and
$\int_{0}^{1}u^2|\hat{\psi}(u)| du \leq
\frac{C(\phi)\log(eT)}{T^2}$.
Also, on making the change of variable $uT \mapsto u$ in the first two
integrals in \eqref{prop22} and recalling the definition of $A(u)$ we
immediately see that these integrals are, respectively, the same as
the second and third integrals on the right hand side of
\eqref{prop2}. Since $C_3(\phi) \leq C(\phi)$, the preceding remarks
together with \eqref{prop24} and \eqref{prop23} gives \eqref{prop2}.

 Let us now simplify \eqref{prop2} further.
 Note that we  have $|A(u)| \leq 2 B(\kappa)x^{\kappa}e^{\kappa}$ when $|u| \leq 1$ by \eqref{cah}. The triangle inequality gives 
\begin{equation}
\label{cor11}
\frac{\kappa}{T} \left| \int_{0}^{T} \left(A(\frac{u}{T})  -A(-\frac{u}{T})\right) u\hat{\phi}(u) \; du \right| \leq \frac{ 4 \kappa B(\kappa)x^{\kappa}e^{\kappa}}{T}\int_{0}^{T} u |\hat{\phi}(u)| \leq \frac{ 8C(\phi) \kappa B(\kappa)x^{\kappa}e^{\kappa}}{T},
\end{equation} 
since $\int_{0}^{T} u |\hat{\phi}(u)|\leq C_0(\phi) +C_3(\phi)$. By
the definition of $A(u)$, the integrand in the first term of
\eqref{cor11} is the same as that in the second integral on the right
hand side of \eqref{prop2}. Thus the corollary follows from the above
estimate and \eqref{prop2}, on noting
that$\frac{\log(eT)}{T^2} \leq \frac{1}{T}$ when $T \geq 1$.
\end{proof}

\begin{remark} In basic applications it is useful to further simplify
  the second term on the right hand side of
  \eqref{cor1}. Thus suppose that $\phi, T$ satisfy the conditions of
  the above corollary with $m = n+1$, $n \geq 1$ and let us for
  brevity set
  $ E(\frac{u}{T}) = \sum_{xe^{-\frac{u}{T}} < n \leq
    xe^{\frac{u}{T}}} a_n {\rm sgn}(x-n)$.
  Then on rewriting $ E(\frac{u}{T})$ as
  $2A(0) - A(\frac{u}{T}) - A(-\frac{u}{T})$ and using the Cahen bound
  \eqref{cah} as above we get
\begin{equation}
\label{ex1}
\int_{T^{\frac{1}{n}}}^{T}  E(\frac{u}{T}) \hat{\phi}(u) \; du \ll_{\kappa, \phi}  \frac{x^{\kappa}}{T},
\end{equation}
since $\int_{T^{\frac{1}{n}}}^{T} u |\hat{\phi}(u)| \; du \ll_{\phi} {1}/{T}$. Also, by the triangle inequality we have 
 \begin{equation}
 \label{ex2}
 \int_{0}^{T^{\frac{1}{n}}}  E\bigl(\frac{u}{T}\bigr) \leq 2C(\phi)
    \max_{0\le \xi\le eT^{\frac{1}{n}-1}}
    \biggl(
    \biggl|\sum_{1< \frac{n}{x}\le 1+\xi}
    a_n
    \biggr|
    +
    \biggl|\sum_{1\leq \frac{x}{n} < 1+\xi}
    a_n
    \biggr| 
    \biggr),
  \end{equation}
  since $1 +\xi = e^{\frac{u}{T}}$ implies $\xi \leq \frac{eu}{T}$
  when $0 \leq u \leq T$, by the mean value theorem, and we have
  $\int_{0}^{T} |\hat{\phi}(u)|\leq C_0(\phi) +C_2(\phi)$.  It follows
  from \eqref{ex1} and \eqref{ex2} that the sum of the second and
  third terms on the right hand side of \eqref{cor1} can be replaced
  with
\begin{equation}
2C(\phi)
    \max_{0\le \xi\le eT^{\frac{1}{n}-1}}
    \biggl(
    \biggl|\sum_{1\le \frac{n}{x}\le 1+\xi}
    a_n
    \biggr|
    +
    \biggl|\sum_{1\le \frac{x}{n}\le 1+\xi}
    a_n
    \biggr|
    \biggr) + {\mathcal O}_{\kappa}\left(\frac{x^{\kappa}}{T}\right).
\end{equation} 
\end{remark}

When used with a suitable $\phi$, for instance with
$\phi = {\mathfrak p}_{3}(t;1)$ of \eqref{exp3}, Corollary~\ref{pp2}
is of similar strength to Theorem 1 of Wolke \cite{wol}. The presence
of the kernel $\phi$ dispenses with the delicate analysis required for
the proof of Theorem 2 of \cite{wol}. Also, Corollary \ref{pp2} merits
comparison with Theorem 2.1 of Liu \& Ye \cite{ly}.  In addition to
the facts that \eqref{prop2} has sums of the $a_n$ rather than $|a_n|$
and is valid for $\kappa > \max(0,\sigma_{c})$, we note that the error
term in \eqref{prop2} has a ${1}/{T}$ rather than essentially
$1/\sqrt{T}$ in Theorem 2.1 of \cite{ly}.

It is perhaps pertinent here to remark that there is a small mistake
in Theorem~1 of \cite{wol}: in inequality (2.5) therein, a factor
$(\frac{x}{n})^{\sigma}$ appears to be missing. This has the
consequence that Theorem 2 of \cite{wol} is valid only for
$T \geq \log x$, a restriction that is of no consequence for the
applications. Theorem 1 of \cite{kp} must therefore also be read with
the same restriction (A. Perelli agrees on this point) as it relies on
\cite{wol}.

One may hope to use the symmetry on account of  the factor ${\rm sgn}(x - n)$ in the
first error term of Proposition \ref{pp2}. This is undoubtedly very difficult in general, but see Coppola
and Salerno \cite{cs1} and \cite{cs2} for a treatment. Theorem 1 in Kaczorowski \& Perelli \cite{kp} also gives a formula with a similar symmetry.

\section{Proof of the Theorem}
\label{pap}

With notation as in the statement of Theorem \ref{theo1}, let ${\mathfrak U}_{q}$ be the set of primitive roots modulo $q$, that is, the set of generators of the multiplicative group $({\bf Z}/q{\bf Z})^{*}$, for a prime number $q$ in $[Q, 2Q]$.  If for any integer $n$ we write ${\tilde n}$ to denote the image of $n$ modulo $q$ and write ${1}_{{\mathfrak U}_q}$ for the characteristic function of ${\mathfrak U}_{q}$, then we have that  
\begin{equation}
\label{one}
{1}_{{\mathfrak U}_q}({\tilde n}) = \sum_{\chi\, {\rm mod}\, q} c_{q}(\chi) \chi(n) \, ,
\end{equation} 
for all integers $n$, where the sum runs over all Dirichlet characters $\chi$ modulo $q$ with $c_{q}(\chi)$  defined to be $\frac{1}{\phi(q)}\sum_{a \in {\mathfrak U}_q } \overline{\chi(a)}$. Also, by an application of the Cauchy-Schwarz inequality followed by the Parseval relation for the group $({\bf Z}/q{\bf Z})^{*}$ we get 
\begin{equation}
\label{2}
\sum_{\chi\, {\rm mod}\, q} |c_q(\chi)| \leq \phi(q)^{\frac{1}{2}} \left(\sum_{\chi\, {\rm mod}\, q} |c(\chi)|^2\right)^{\frac{1}{2}} = |{\mathfrak U}_q|^{\frac{1}{2}}.
\end{equation} 
Throughout the remainder of this section ${\bf b} = \{b_n\}_{n \geq 1}$ will denote one of the sequences $\{\Lambda(n)\}_{n \geq 1}$ and $\{\mu(n)\}_{n \geq 1}$. Then for any real number $w \geq 1$ we have 
\begin{equation}
\label{31} 
\sum_{\substack{Q\leq q \leq 2Q, \\ q\, \text{prime}}} \sum_{\substack{1 \leq n \leq w, \\ n \in {\mathfrak U}_q  \, {\rm mod}\, q. }} b_n = \sum_{\substack{Q\leq q \leq 2Q, \\ q\, \text{prime}}} \sum_{\chi\, {\rm mod}\, q} c_q(\chi) \sum_{1 \leq n \leq w} \chi( n) b_n .
\end{equation}
For a given real $x \geq 1$, let us set $y = x^{\theta}$ with
$0 < \theta < 1$. Then on subtracting the contribution from the
principal character $\chi_{0}$ modulo $q$ to the right hand side of
\eqref{31} from both sides of this relation and using the resulting
relation for $w = x$, $w = x+y$ together with triangle inequality we
get
\begin{equation}
\label{4} 
\sum_{\substack{Q\leq q \leq 2Q, \\ q\, \text{prime}}} 
\biggl|
\sum_{\substack{x < n \leq x+y, \\ n \in {\mathfrak U}_q  \, {\rm
      mod}\, q. }} b_n - \frac{|{\mathfrak
    U}_{q}|}{\phi(q)}\sum_{\substack{x < n \leq x+y, \\ (n,q) =1. }}
b_n \biggr|\; 
\leq\; \sum_{\substack{Q\leq q \leq 2Q, \\ q\, \text{prime}}} 
\biggl|\sum_{\substack{\chi\, {\rm mod}\, q, \\ \chi \neq \chi_{0}.}} c_q(\chi) \sum_{x < n \leq x+y} \chi( n) b_n \biggr|.
\end{equation}
We shall presently bound the sum  
\begin{equation}
\label{3}
\Sigma = \Sigma({\bf b}, x,y, Q) = 
\sum_{\substack{Q\leq q \leq 2Q, \\ q\, \text{prime}}} 
\biggr|\sum_{\substack{\chi\, {\rm mod}\, q, \\ \chi \neq \chi_{0}.}} c_{q}(\chi) \sum_{x < n \leq x+y} b_n \chi(n)\biggr|.
\end{equation}
by means of Corollary \ref{pp2}.  To this end, we set
$F(s, \chi) = \sum_{n \geq 1} {b_n \chi(n)}/{n^s}$, which
converges in $\sigma > \frac{1}{2}$ for each $\chi \neq \chi_{0}$
under the GRH for our sequences ${\bf b}$. We then fix a $\epsilon >0$
and set $\kappa = \frac{1}{2} + \epsilon$. Also, we let $\varphi$ be a
positive continuous function supported in $[-1,1]$ and satisfying the
conditions on $\phi$ of Corollary \ref{pp2} with $m =2$. For example
we may take $\phi = {\mathfrak p}_{3}(t;1)$ of \eqref{exp3} . On
applying this proposition we now get
\begin{equation}
 \label{5}
 \begin{split}
 \sum_{1 \leq n \leq w} b_n \chi(n) \; &= \frac{1}{2\pi i} \int_{\kappa-i\infty}^{\kappa + i\infty} F(s, \chi)\varphi\left(\frac{s-\kappa}{2\pi iT }\right) \frac{w^s}{s} \; ds \, \\
            &+\int_{0}^{T} \left(\sum_{we^{-\frac{u}{T}} < n \leq we^{\frac{u}{T}}} b_n\chi(n) {\rm sgn}(w-n)\right) \hat{\varphi}(u) \; du \\
            &+ O\left( \frac{B(\kappa, \chi)w^{\kappa}\log(eT)}{T^2}\right), 
 \end{split}
 \end{equation}
 for all real numbers $w \geq 1$, $T \ge1$ and $\chi \neq
 \chi_{0}$.
 Here
 $B(\kappa, \chi) = \sup_{N \geq 1}|\sum_{1 \leq n \leq N}
 \frac{b_n\chi(n)}{n^{\kappa}}|\ll_{\epsilon} q^{\epsilon}$
 under the GRH for sequences ${\bf b}$ given above, as can be seen by
 integrating by parts using Theorem 15.5 of \cite{iwa}. We now note
 the following lemma, which allows us to take advantage of the
 cancellation in the sums on the right hand side of the above
 relation.

\begin{lemma} \label{lem1} Let $w \geq 1$ be a real number. Then if $u_n = b_n$ or $u_n =  b_n {\rm sgn}(w-n)$ for all $n \geq1$,  with $\{b_n\}_{n \geq 1}$ as above, we have 
\begin{multline}
\label{8}
\int_{0}^{T} \sum_{\substack{Q\leq q \leq 2Q, \\ q\, \text{prime}}} 
\biggl|\sum_{\substack{\chi\, {\rm mod}\, q, \\ \chi \neq \chi_{0}.}}
c_{q}(\chi) \sum_{we^{-\frac{u}{T}} < n \leq we^{\frac{u}{T}}}
u_n\chi(n) \hat{\varphi}(u)\biggr|
 \; du  
 \\\ll \max_{1 \leq n \leq ew}|u_n| \Bigl(\frac{w\log(eT)}{T} + Q\Bigr) Q^{\frac{1}{2}},
\end{multline}
where the implied constant depends on $\varphi$ alone.
\end{lemma}

\begin{proof}
By means of the triangle inequality and the Cauchy-Schwarz inequality  we have 
\begin{equation}
\label{6}
\begin{split}
&\Biggl(\int_{0}^{T} \sum_{\substack{Q\leq q \leq 2Q, \\ q\, \text{prime}}} |\sum_{\substack{\chi\, {\rm mod}\, q, \\ \chi \neq \chi_{0}.}} c_{q}(\chi) \sum_{we^{-\frac{u}{T}} < n \leq we^{\frac{u}{T}}} u_n\chi(n) \hat{\varphi}(u)| \; du \Biggr)^2 \\
&\leq  
\biggl(\sum_{\substack{Q\leq q \leq 2Q, \\ q\, \text{prime}}} \sum_{\substack{\chi\, {\rm mod}\, q, \\ \chi \neq \chi_{0}.}} |c_{q}(\chi)|^2 \int_{0}^{T} | \hat{\varphi}(u)| \; du \biggr) \times \\
&\qquad\biggl(\int_{0}^{T} \sum_{\substack{Q\leq q \leq 2Q, \\ q\,
    \text{prime}}} \sum_{\substack{\chi\, {\rm mod}\, q, \\ \chi \neq
    \chi_{0}.}}  
\Bigl| \sum_{we^{-\frac{u}{T}} < n \leq we^{\frac{u}{T}}} u_n \chi(n)\Bigr|^2\, |\hat{\varphi}(u)| \; du \biggr).
\end{split}
\end{equation}
Since $\hat{\varphi}$ is integrable on ${\bf R}$ and $\sum_{\substack{\chi\, {\rm mod}\, q, \\ \chi \neq \chi_{0}.}} |c_{q}(\chi)|^2 \leq {|{\mathfrak U}_{q}| }/{\phi(q)}$ by the Parseval relation, the first of the two bracketed expressions on the right hand side of the above relation does not exceed $\|\hat{\varphi}\|_{1} \sum_{Q \leq q \leq 2Q}\frac{|{\mathfrak U}_{q}|}{\phi(q)} \ll Q$. We estimate the second expression using a variant of the large sieve inequality for characters. Indeed, when  $w\ge 1$ and $u>0$, the number of integers in $(we^{-\frac{u}{T}},we^{\frac{u}{T}}]$ is at most 
$\frac{2wue^{\frac{u}{T}}}{T}+1$. Then it follows from this inequality that 
\begin{align}
  \sum_{\substack{Q\le q\le 2Q,\\\text{$q$ prime}}}
  \sum_{\substack{\chi \mod q, \\ \chi\neq \chi_0}}
  &\biggl|\sum_{we^{-\frac{u}{T}} < n \leq we^{\frac{u}{T}}}u_n\chi(n)\biggr|^2
  \notag
  \\&\le
  \max_{we^{-\frac{u}{T}} < n \leq we^{\frac{u}{T}}}|u_n|^2
  \Bigl(\frac{2wue^{\frac{u}{T}}}{T}+1\Bigr)\Bigl(\frac{2wue^{\frac{u}{T}}}{T}+4Q^2\Bigr)\notag
  \\&\label{bound1}
  \ll \max_{we^{-\frac{u}{T}} < n \leq we^{\frac{u}{T}}}|u_n|^2\Bigl(\frac{2wue^{\frac{u}{T}}}{T}+Q\Bigr)^2.
\end{align}     
On noting that $\int_{0}^{T}u^{2} |\hat{\varphi}(u)| \; du \leq C(\phi) \log(eT)$ and  $\int_{0}^{T}|\hat{\varphi}(u)| \; du \leq 2C(\phi)$ we conclude that the second expression in the brackets on the right hand side of \eqref{6} is majorised by 
\begin{equation*}
 \max_{1 \leq n \leq ex}|u_n|^2 \int_{0}^{T} \Bigl(\frac{w^2u^2e^{\frac{2u}{T}}}{T^2}+Q^2\Bigr) |\hat{\varphi}(u)| \; du  \ll  \max_{1 \leq n \leq ew}|u_n|^2 \Bigl(\frac{w^2\log(eT)}{T^2} + Q^2\Bigr)\;.
\end{equation*}
The lemma now follows on substituting the preceding bounds into \eqref{6} and passing to square roots.
\end{proof}

We sum the absolute values of both sides of \eqref{5} over the characters $\chi \neq \chi_{0}$ and the primes $q$ in $[Q, 2Q]$. We then estimate the second and third terms on the right hand side of the resulting relation using Lemma \ref{lem1}. On using \eqref{2} to bound the error term of this relation we conclude that    
\begin{equation}
 \label{81}
 \begin{split}
\sum_{\substack{Q\leq q \leq 2Q, \\ q\, \text{prime}}} \sum_{\substack{\chi\, {\rm mod}\, q, \\ \chi \neq \chi_{0}.}} \sum_{1 \leq n \leq w} b_n \chi(n) \; &= \frac{1}{2\pi i} \int_{\kappa-i\infty}^{\kappa + i\infty} \sum_{\substack{Q\leq q \leq 2Q, \\ q\, \text{prime}}} \sum_{\substack{\chi\, {\rm mod}\, q, \\ \chi \neq \chi_{0}.}}F(s, \chi)\varphi\left(\frac{s-\kappa}{2\pi iT }\right) \frac{w^s}{s} \; ds \, \\
            &+ O\left( \max_{1 \leq n \leq ew}|b_n| \Bigl(\frac{w\log(eT)}{T} + Q\Bigr) Q^{\frac{1}{2}}\right)\\
            &+ O_{\epsilon}\left( \frac{Q^{\epsilon}\log(eT)w^{\frac{1}{2}+\epsilon}\sum_{Q \leq q \leq 2Q} |{\mathfrak U}_q|^{\frac{1}{2}}}{T^2}\right) 
 \end{split}
 \end{equation}
for all real $w \geq 1$. We apply this with $w = x$ and $w = x+y$, subtract and recall the definition of $\Sigma$ to obtain by means of the triangle inequality that 
\begin{equation}
 \label{9}
\begin{split}
\Sigma \; &\leq 
\frac{1}{2\pi } \int_{\kappa-i\infty}^{\kappa + i\infty} \sum_{\substack{Q\leq q \leq 2Q, \\ q\, \text{prime}}} \sum_{\substack{\chi\, {\rm mod}\, q, \\ \chi \neq \chi_{0}.}}|c_{q}(\chi)||F(s, \chi)| \biggl|\varphi\left(\frac{s-\kappa}{2\pi iT }\right)\biggr| \left|\frac{(x+y)^s- x^s}{s}\right| \; dt \, \\
            &+ O\left( \max_{1 \leq n \leq 2ex}|b_n| \Bigl(\frac{x\log(eT)}{T} + Q\Bigr) Q^{\frac{1}{2}}\right)
            + O_{\epsilon}\left( \frac{Q^{\frac{3}{2} +\epsilon}x^{\frac{1}{2}+\epsilon}\log(eT)}{T^2}\right). 
\end{split}
\end{equation}
On the GRH we have the classical Lindel{\"o}f bound $|F(s,\chi)| \ll_{\epsilon} (q+q|t|)^{\epsilon}$, by \cite{iwa}, Theorem 5.17 and Corollary 5.19.
Also, for $s = \kappa +it$ we have $|\frac{(x+y)^s- x^s}{s}| \leq \min(\frac{3x^{\kappa}}{|s|}, x^{\kappa-1}y)$ by a trivial estimate and the mean value theorem. Further, $\varphi(\frac{s-\kappa}{2\pi iT}) =0$ when $|t| \geq T$. On combining these remarks with \eqref{2} and assuming that $T \leq x$, we see that  
\begin{equation}
\label{10}
\begin{split}
 &\int_{\kappa-i\infty}^{\kappa + i\infty} \sum_{\substack{Q\leq q \leq 2Q, \\ q\, \text{prime}}} \sum_{\substack{\chi\, {\rm mod}\, q, \\ \chi \neq \chi_{0}.}}|c_{q}(\chi)||F(s, \chi)| |\varphi\left(\frac{s-\kappa}{2\pi iT }\right)| \left|\frac{(x+y)^s- x^s}{s}\right| \; dt \,\\
 &\ll_{\epsilon} Q^{\frac{3}{2} +\epsilon} x^{\epsilon} \int_{-T}^{T} \min(\frac{3x^{\kappa}}{|\kappa+it|}, x^{\kappa-1}y) \, dt \ll_{\epsilon}  \frac{y Q^{\frac{3}{2}}}{x^{\frac{1}{2}}}\min(T, \frac{x}{y}) (xQ)^{\epsilon} 
\end{split}
\end{equation}
Using this in \eqref{8} and noting that $\log(eT) \ll_{\epsilon} x^{\epsilon}$ we finally obtain 
    \begin{equation}
    \label{11}
   \begin{split}
   \Sigma \; &\ll_{\epsilon} (xQ)^{\epsilon}Q^{\frac{3}{2}}\left(\frac{y }{x^{\frac{1}{2}}}\min(T, \frac{x}{y}) + \frac{x}{TQ} + 1  + \frac{x^{\frac{1}{2}}}{T^2}\right), 
   \end{split}
   \end{equation}
since for our choices of the sequence ${\bf b}$ we certainly have $\max_{1 \leq n \leq 2ex}|b_n| \ll_{\epsilon} x^{\epsilon}$. We set $T=\frac{x^{3/4}}{Q^{1/2}y^{1/2}}$ and  note that $\min(T,x/y)\le T$, which
holds since $y\le Qx^{\frac{1}{2}}$. Then on combining \eqref{11} with \eqref{4} and \eqref{3} we get the bound
\begin{equation}
\label{12}
\sum_{\substack{Q\leq q \leq 2Q, \\ q\, \text{prime}}} |\sum_{\substack{x < n \leq x+y, \\ n \in {\mathfrak U}_q  \, {\rm mod}\, q. }} b_n - \frac{|{\mathfrak U}_{q}|}{\phi(q)}\sum_{\substack{x < n \leq x+y, \\ (n,q) =1. }} b_n |
   \ll_{\epsilon} yQ\, (Qx)^{\epsilon}
  \biggl(\frac{x^{\frac{1}{4}}}{y^{\frac{1}{2}}}
  +\frac{Q^{\frac{1}{2}}}{y}
  +\frac{Q^{\frac{3}{2}}}{ x}\biggr).
\end{equation} 

\subsection{The case of the primes}

We now take $b_n = \Lambda(n)$ in \eqref{12} and verify the first conclusion of Theorem \ref{theo1}. In effect, since $Q \leq x$, in this case \eqref{12} can be rewritten as 
\begin{equation*}
  \frac{\log Q}{Q}\sum_{\substack{Q\le q\le 2Q,\\\text{$q$ prime}}}
  \biggl|\sum_{\substack{x < n \leq x+y, \\ n \in {\mathfrak U}_q  \, {\rm mod}\, q.}}\Lambda(n)
  -
  \frac{|\mathfrak{A}_q|}{\phi(q)}
  \sum_{\substack{x < n \leq x+y, \\ (n,q) =1. }}\Lambda(n)
  \biggr|
 \ll_{\epsilon}
  y\,(Qx)^{\epsilon}
  \biggl(\frac{x^{1/4}}{y^{1/2}}
  +\frac{Q^{1/2}}{y}
  +\frac{Q^{3/2}}{ x}\biggr).
\end{equation*}
Under the RH we have $\sum_{\substack{x < n \leq x+y}}\Lambda(n) = y +O(x^{\frac{1}{2}}(\log x)^2))$. The trivial estimate  for the contribution from $n = p^k$, with $p$ prime and $k \geq 2$, to the sums inside the absolute value on the left hand side is $O( x^{\frac{1}{2}}(\log x)^2)$. Since $Q \leq x$, the condition $(n,q) =1$ on the left hand side can be dropped when $n$ is a prime. These remarks yield 
 \begin{equation*}
   \frac{\log Q}{Q}\sum_{\substack{Q\le q\le 2Q,\\\text{$q$ prime.}}}
   \biggl|\sum_{\substack{x < p \leq x+y, \\ p \in {\mathfrak U}_q  \, {\rm mod}\, q,\\ \text{$p$ prime.}}} \log p
   -
   \frac{|\mathfrak{A}_q|y}{\phi(q)}
   \biggr|
  \ll_{\epsilon}
   y\,(Qx)^{\epsilon}
   \biggl(\frac{x^{1/4}}{y^{1/2}}
   +\frac{Q^{1/2}}{y}
   +\frac{Q^{3/2}}{ x}\biggr) + x^{\frac{1}{2}}(\log x)^2.
 \end{equation*}
 With $y = x^{\frac{1}{2} +\delta}$ we have $x^{\frac{1}{2}}(\log x)^2 \ll yx^{\epsilon}x^{-\frac{\delta}{2}}$. Since also $Q\le x^{\frac{2}{3+2\theta}}$ we then get 
  \begin{equation}
  \label{16}
    \frac{\log Q}{Q}\sum_{\substack{Q\le q\le 2Q,\\\text{$q$ prime.}}}
    \biggl|\sum_{\substack{x < p \leq x+y, \\ p \in {\mathfrak U}_q  \, {\rm mod}\, q,\\ \text{$p$ prime.}}} \log p
    -
    \frac{|\mathfrak{A}_q|y}{\phi(q)}
    \biggr|
   \ll_{\epsilon} yx^{\epsilon}\Bigl(x^{-\delta/2}
     +x^{\frac{1}{3+2\theta}-\frac12-\delta}
     +x^{\frac{3}{3+2\theta}-1}\Bigr).
      \end{equation}
We set $\eta=\min(\frac16,\frac\delta2,\frac{2\theta}{3})$ and choose
$\epsilon \le \frac{\eta}{3}$ to find that
\begin{equation*}
  x^{\epsilon}\Bigl(x^{-\delta/2}
  +x^{\frac{1}{3+2\theta}-\frac12-\delta}
  +x^{\frac{3}{3+2\theta}-1}\Bigr)
  \ll x^{\frac{-2\eta}{3}}.
\end{equation*} 
For any prime number $q$ we have $|\mathfrak{A}_q| = \phi(p-1)$ and ${\phi(q)= p-1}$. Thus if ${\mathcal S}$ is the set of primes $q$ in $[Q, 2Q]$ such that 
\begin{equation}
\label{18}
  \biggl|\sum_{\substack{x-y\le p\le x,\\ p\in\mathfrak{A}_q}}\log p
  -
  \frac{\varphi(q-1)y}{q-1}
  \biggr|\ge \frac{\varphi(q-1)y}{q-1}
  x^{-\eta/3}
\end{equation}
then it follows from \eqref{16} and $\frac{\varphi(q-1)y}{q-1} \gg \frac{1}{\log\log Q}$ when $q \geq Q$ that  
\begin{equation}
\label{19}
 \frac{ |\mathcal{S}|}{Q} \ll x^{-\eta/3}\ll Q^{-\eta/2} \leq Q^{\frac{1}{12}}
  ,
\end{equation} 
which yields the desired conclusion of the theorem after removing the weights $\log p$ in the usual fashion.   

\subsection{The case of the M{\"o}bius function}

Here we set $b_n = \mu(n)$ in \eqref{12} and carry out the details just as in the preceding case, taking note of the simplification afforded by the fact that in this case there is  no main term and no prime powers in the support of the function $\mu$.

\section{Moving the Line of Integration}
\label{mli}

Our first purpose here is to record the proposition below that describes the effect of
``moving the line of integration" in the integrals over the line $\sigma = \kappa$ on the right hand sides of \eqref{prop2} and \eqref{cor1} when  $\phi$ is a given continuous positive compactly supported piecewise polynomial function. 

\vspace{2mm}
\noindent
It will be convenient here to use both $s = \sigma + it$ and
$z = u +iv$ to denote complex numbers. Also, we shall suppose that the
support of $\phi$ is in $[-U, U]$ for some $U >0$. Further, let
$-U= u_{1} < u_{2} < \ldots < u_{m} =U$ be such that the restriction
of $\phi$ to the real interval $[u_j, u_{j+1})$ agrees with that of a
polynomial $\widetilde{\phi}_{j}$ defined on ${\bf C}$, for
$1 \leq j \leq m-1$. We will assume that $u_j \neq 0$ and let
$\eta \leq |u_j|$ for all $j$.  Let $V >0$ be a positive real number
and let $M$ satisfy $|\widetilde{\phi}^{\prime}_{j}(z)| \leq M$ for
all $z = u +iv$ in the rectangle $-U \leq u \leq U$ and
$0 \leq v \leq V$ and $1 \leq j \leq m-1$. Finally, we define
$\widetilde{\phi}(z)$ for $z = u+iv$ with $u \in [-U, U)$ by
$\widetilde{\phi}(z) = \widetilde{\phi}_{j}(z)$ where $j$ is the
unique index such that $u\in [u_j, u_{j+1})$ .

\begin{proposition}
  \label{pp3} With notation as above, let $ \kappa^{\prime}, \kappa$
  be such that $ 0 < \kappa - \kappa^{\prime} \leq 2 \pi V$. Also, let
  $F$ be a meromorphic function on a neighbourhood of the closed
  rectangle with vertices $\kappa^{\prime} \pm 2\pi iUT$ and
  $\kappa \pm 2\pi iUT$, for some $T \geq 1 $. Suppose further that if
  ${\mathcal A}$ is the set of poles of $s \mapsto {F(s)}/{s}$ in
  this neighbourhood then ${\rm Re}(a) \neq \kappa, \kappa^{\prime}$
  and ${\rm Im}(a) \neq 2\pi u_j$ for all $a$ in ${\mathcal A}$ and
  $1 \leq j \leq m$. Then we have that
\begin{multline}
\label{prop31}
\frac{1}{2\pi i} \int_{\kappa-i\infty}^{\kappa + i\infty} F(s)\phi\left(\frac{s-\kappa}{2\pi iT }\right) \frac{x^s}{s} \; ds \, 
= \frac{1}{2\pi i} \int_{\kappa^{\prime}-i\infty}^{\kappa^{\prime} +
  i\infty} F(s)\widetilde{\phi}\left(\frac{s-\kappa}{2\pi iT }\right)
\frac{x^s}{s} \; ds \, 
\\+ \,
\sum_{a \in {\mathcal A}} \,{\rm
  Res}\left(\frac{\widetilde{\phi}(s)F(s)x^s}{s}\right)_{s =a} 
 \\+ \mathcal{O}^{*}\left( \frac{M(\kappa- \kappa^{\prime})x^{\kappa}}{4 \pi^3 \eta T^2} \sum_{1 \leq j \leq m}\int_{\kappa}^{\kappa^{\prime}} |F(\sigma +2\pi iu_jT)| d\sigma \right).
\end{multline}
\end{proposition}

\begin{proof}
  For $1 \leq j \leq m-1$, the function $G$ with
  $G(s) = \frac{F(s)\widetilde{\phi}_{j}\left(\frac{s-\kappa}{2\pi iT
      }\right) x^s}{s}$
  is meromorphic in a neighbourhood of the closed rectangle
  ${\mathcal R}_{j}$ with vertices $\kappa^{\prime} + 2\pi iu_{j+1}T$,
  $\kappa^{\prime} + 2\pi iu_{j}T$ and $\kappa + 2\pi iu_{j}T$,
  $\kappa + 2\pi iu_{j+1}T$, with no poles on the boundary of this
  rectangle. On applying the residue theorem to $G$ on each
  ${\mathcal R}_j$ oriented anticlockwise for $1 \leq j \leq m-1$ and
  adding the resulting relations, we see that \eqref{prop31} follows
  if we show that the sum of the integrals of $G$ along the oriented
  horizontal sides of the ${\mathcal R}_j$ is majorised by the error
  term in \eqref{prop31}. This reduces to verifying for
  $1 \leq j \leq m$ the inequality
\begin{equation}
\label{prop32}
\left|\frac{1}{2\pi i} \int_{\kappa^{\prime}+ 2\pi i u_j T}^{\kappa^{\prime} + 2 \pi i u_j T} F(s)x^s P_j\left(\frac{s-\kappa}{2\pi iT }\right) \frac{ ds}{s} \right|\, \leq \frac{M(\kappa- \kappa^{\prime})x^{\kappa}}{4 \pi^2 \eta T^2} \int_{\kappa}^{\kappa^{\prime}} |F(\sigma +2\pi iu_jT)| d\sigma,
\end{equation}
where $P_1 (s) = \widetilde{\phi}_{1}(s)$,
$P_m(s) = \widetilde{\phi}_{m-1}(s)$ and
$P_j(s) = \widetilde{\phi}_{j}(s)-\widetilde{\phi}_{j-1}(s)$ for
$2 \leq j \leq m-1$. For $1 \leq j \leq m$ and
$\sigma \in [\kappa^{\prime}, \kappa]$, let us set
$a_j(\sigma) = P_j\left(\frac{\sigma + 2\pi i u_jT-\kappa}{2\pi iT
  }\right)$.
Then since $\phi$ is continuous and supported in $[-U,U]$, we have
$a_j(\kappa) = P_j(u_j) =0$ for each $j$. Thus the mean value theorem
applied to $\sigma \mapsto a_{j}(\sigma)$, which is a continuously
differentiable function, gives
\begin{equation}
\label{prop33}
\biggl|P_j\left(\frac{\sigma + 2\pi i u_jT-\kappa}{2\pi iT }\right)\biggr| = |a_{j}(\sigma) -a_j(\kappa)| \leq \frac{2M(\kappa - \kappa^{\prime})}{2 \pi T}, 
\end{equation} 
for all $\sigma \in [\kappa^{\prime}, \kappa]$ and $1 \leq j \leq
m$.
Here we have used
$0 \leq \frac{\kappa - \kappa^{\prime}}{2 \pi T} \leq V$, since
$T \geq 1$. Also, for $\sigma \in [\kappa^{\prime}, \kappa]$ we have
$|x^{\sigma+ 2 \pi i u_jT}| \leq x^{\kappa}$ and
$|\sigma+ 2 \pi i u_jT| \geq 2 \pi \eta T$ for $1 \leq j \leq
m$.
These bounds together with an application of the triangle inequality
to the left hand side of \eqref{prop32} verify this inequality.
\end{proof}

We now describe a convenient family test functions that may be used
for $\phi$ in our formulae. Let us we set $\delta >0$ and $m \geq 1$,
an integer. Also, we will write $1_{[a,b]}$ for the characteristic
function of the interval $[a,b]$ and $1_{[a,b]}^{(*m)}$ for the $m$-th
convolution of $1_{[a,b]}$ with iteslf.  Then we define even function
${\mathfrak p}_{m}(t;\delta)$ by
\begin{equation}
  \label{defpmfm}
  {\mathfrak p}_{m}(t;\delta)=\int_{-\infty}^{\infty}
  1_{[-1,1]} \biggl(\frac{u}{1+\frac{1}{2}\delta}\biggr)
  1_{[-1,1]}^{(*m)} \biggr(\frac{t-u}{\frac{1}{2}\delta/m}\biggr)
  \frac{m\,du}{2^m \frac{1}{2}\delta}.
\end{equation}
It is easily seen that ${\mathfrak p}_{m}(t;\delta)$ is piecewise
polynomial of class $C^{m}$ and that its support lies in
$[-(1+\delta),1+\delta]$. Moreover, we have
${\mathfrak p}_{m}(t;\delta)=1$ when $|t|\le 1$ and
$ 0\le {\mathfrak p}_{m}(t;\delta)\le 1$ for all real $t$. Finally, it
immediately follows from basic properties of the Fourier transform
that
\begin{equation}
  \label{pmffourier}
  \hat{\mathfrak p}_{m}(u;\delta)
  =\frac{\sin(\pi(2+\delta)u)}{\pi u}
  \biggl(\frac{m\sin(\pi \delta u/m)}{\pi \delta u}\biggr)^m
\end{equation}
and that $\hat{\mathfrak p}_{m}(-u;\delta)=\hat{\mathfrak p}_{m}(u;\delta)$. We end this note by explicitly describing ${\mathfrak p}_{3}(t;1)$ :
\begin{equation}
  \label{exp3}
  \displaystyle\quad
  {\mathfrak p}_{3}(t;1)=
  \begin{cases}
    0&\text{when $|t|\ge 3/2$},
    \\\displaystyle
    3(3-2t)^3/16
    &\text{when $7/6\le |t|\le 3/2$},
    \\\displaystyle
    (36t^3-108t^2+99t-25)/4&
    \text{when $5/6\le |t|\le 7/6$},\\
    (-9t^3+27t^2-27t+25)/16&\text{when $1/2\le |t|\le 5/6$},\\
    1&\text{when $|t|\le 1/2$}.
  \end{cases}
  \quad
\end{equation}

\vspace{3mm}
\noindent
{\bf Acknowledgement :}
The authors are grateful to Herv\'e Queffelec for the many interesting
discussions on Fourier and Mellin transforms and Perron's
formula. Thanks are also due to Eero Saskman for sharing his
ideas. This work was put in final form with support from the CEFIPRA
project 5401-1. The first author gratefully acknowledges the
facilities provided to him by the Universit{\'e} Aix-Marseille during
his visit under the aegis of the said project.

\end{document}